\documentclass[prl,showpacs]{revtex4-1}
\usepackage{graphicx, amsmath, amssymb}

\newtheorem{theorem}{Theorem}

\newtheorem{algorithm}{Algorithm}
\newtheorem{definition}{Definition}
\newtheorem{example}{Example}
\newtheorem{remark}{Remark}

\begin{document}

\title{Characterizing the Topography of Multi-dimensional Energy Landscapes}
\author{H.\,Lydia Deng}
\address{Landmark Graphics Corp \\
Highlands Ranch, Colorado USA}

\author{John A.\,Scales}
\address{Department of Physics \\
Colorado School of Mines \\
Golden, Colorado 80401, USA}


\begin{abstract}
  A basic issue in optimization, inverse theory,
  neural networks, computational chemistry and many other problems is
  the geometrical characterization of high dimensional functions.  
  In inverse
  calculations one aims to characterize the set of models that fit the
  data (among other constraints).  If the data misfit function is
  unimodal then one can find its peak by local optimization methods
  and characterize its width (related to the range of data-fitting
  models) by estimating derivatives at this peak.  On the other hand,
  if there are local extrema, then a number of interesting and
  difficult problems arise.  Are the local extrema important compared
  to the global or can they be eliminated (e.g., by smoothing) without
  significant loss of information?  Is there a sufficiently small
  number of local extrema that they can be enumerated via local
  optimization?  What are the basins of attraction of these local
  extrama?  Can two extrema be joined by a path that never goes
  uphill?  Can the whole problem be reduced to one of enumerating the
  local extrema and their basins of attraction?  For locally
  ill-conditioned functions, premature convergence of local
  optimization can be confused with the presense of local
  extrema. Addressing any of these issues requires topographic
  information about the functions under study.  But in many
  applications these functions may have hundreds or thousands of
  variables and can only be evaluated pointwise (by some numerical
  method for instance).  In this paper we describe systematic 
  (but generic) methods
  of analysing the topography of high dimensional functions using local
  optimization methods applied to randomly chosen starting models.  We
  provide a number of quantitative measures of function topography
  that have proven to be useful in practical problems along with error
  estimates.
\end{abstract}

\maketitle

\section{What Makes an Optimization Problem Hard?}

We consider the problem of optimizing a function 
$F$ (the {\it objective} or {\it cost} function) mapping 
${\cal M} \subset {\bf R}^N$ into ${\cal Y} \subset {\bf R}$.
We refer to ${\cal M}$ as the {\it model  space}, and each point in
the model space, ${\bf m}$, is a {\it model}. Depending on the
application, the goal may be to find the global extremum of  $F$, a
single local extremum, or a collection of local extrema. 
In this paper we will assume that optimization refers to minimization,
whether local or global.

There is no generally agreed upon characterization of what makes an
optimization problem hard.  Hardness has to do partly with our goals
--- do we need a global extremum or will a good local extremum do;
partly with the structure of the function --- does it have lots of
local extrema, how broad are the basins associated with these extrema;
and partly with the dimensionality of the problem --- exhaustive
search will be infeasible except for low-dimensional problems.

In many applications, however, the function $F$ cannot be expressed in
closed form in terms of elementary functions, but can only be
evaluated point-wise by computer programs. Such problems arise in many
fields. Some of the most widely studied include the {\it spin-glass}
problem, the {\it traveling-salesman} problem (TSP), and the  {\it
residual statics} problem of exploration seismology. 

\subsection{Global Search Strategies}

If the structure of function $F$ is unknown, optimization is
fundamentally a matter of search in the model space. In order to be 
able to treat such a broad variety of situations, we begin with an
abstract statement of a search algorithm. Here, we use the notation
${\bf {\vec m}}^t$ to represent a population of candidate models at
the time step $t$.
\begin{algorithm}{\bf General Search (GS)} 
 ${\bf {\vec m}}={\mathit GS}(F,P,{\bf T},S)$ \\
Let $F:{\cal M}\subset{\bf R}^N \rightarrow {\cal Y}\subset{\bf R}$,
 $P \equiv {\bf {\vec m}}^0 = \{{\bf m}^0_k\}_{k=1,\cdots,K}$ be an
initial population of models, where ${\bf m}^0_k \in {\cal M}$ and 
$K \ge 1$, ${\bf T}$  a {\it transition operator}, and  $S$ a stopping
criterion.
\begin{enumerate}
\item Iteratively apply the transition operator to generate a new
  population of models at each iteration, so that
  ${\bf {\vec m}}^t = {\bf T} \, {\bf {\vec m}}^{t-1};$ 
\item Repeat (1) until $S$ is satisfied. The final set of models 
 $ {\bf {\vec m}}$ are the output of the search.
\end{enumerate}
\label{al:gs}
\end{algorithm}

Any searching process can be considered as an evolution of a
population of models (possibly a single model) in the $N$-dimensional 
model space. The transition operator  ${\bf T}$ is the rule that
determines to which models the population evolves from the previous
population. Here, we assume that the transition operator  ${\bf T}$ is
independent of the time step $t$, which is the case in most
algorithms. Different optimization algorithms differ by the strategies
in choosing the initial population $P$ and the rules of transition
from one population of models to another, ${\bf T}$.

Among the searching methods defined via Algorithm~\ref{al:gs}, there
are two extreme strategies, {\it hill-climbing} (HC) and {\it uniform
Monte Carlo} (UMC). HC search is a local descent search applied to a
single model (population size $K=1$). An initial model  $P={\bf m}^0$
is selected (possibly at random) and the transition operators ${\bf T}
= {\bf T}_{local}$ are deterministic operators, such as conjugate
gradient, quasi-Newton,  or downhill simplex, which follow a path
downhill as far as possible. For objective functions containing more
than one local extremum ({\it multi-modal}), the result of HC strongly
depends on the choice of the initial model ${\bf m}^0$. UMC, on the
other  hand, selects points with uniform probability in the model
space. The transition operation  ${\bf T}$ is simply the selection of
new points at random and therefore makes no use of information from
previous generations. Thus, if there are $N$ parameters and each of
them can take $m$ possible values, the probability of finding a
particular model is proportional to $m^{-N}$ for each function
evaluation using UMC.

Search strategies have been developed that yield a compromise between
these two extremes; almost all of these incorporates stochastic
elements, especially in the construction of transition operators. It
is important for the success of global searches that the transition 
operators make the best use of information provided by the current
samples while avoiding being trapped in local extrema. Among all these
strategies, the most widely used are {\it Simulated Annealing} (SA)
\cite{kirkpatrick}, {\it Genetic Algorithms} (GA) \cite{holland} and
random hill-climbing (RHC), to be defined shortly.

SA and GA searching strategies use stochastic transition operators
 ${\bf T}$ that are biased towards good samples from the previous
generations. Many variations of SA and GA can be found in the
literature \cite{aarts:book,vanlaarhoven,goldberg:book}. Although the
asymptotic convergence results are known for both SA \cite{hajek} and
GA \cite{davis_principe} these results are hardly useful in practice.

RHC searches, on the other hand, apply deterministic transition
operations ${\bf T}$ to a randomly chosen population $P$.
Hence RHC explores locally in multiple areas of objective functions,
and the resulting samples are a set of local/global extrema. This
search algorithm can be described as 
\begin{algorithm}{\bf Random Hill Climbing}
 ${\bf {\vec m}} = {\mathit RHC}(F,K,\epsilon,cmax)$ \\
Let the randomly chosen initial population size be $K$. Let the
stopping criterion $S$ be that either gradients of all samples are
reduced to the tolerance $\epsilon$ or the number of iterations
reaches a maximum  $cmax$. Let ${\bf T}_{local}$ be a local descent
search operator.
\begin{enumerate}
\item Choose initial models 
 $P=\{{\bf m}^0_k\}_{k=1,\cdots,K} \in {\cal M} $ uniformly at random,
where $K \gg 1$;
\item Apply Algorithm~\ref{al:gs}, 
 ${\bf {\vec m}} = {\mathit GS}(F,P,{\bf T}_{local},S)$.
\end{enumerate}
The final population contains $n$ distinct models, 
 ${\bf {\vec m}} = \{{\bf  m}_k\}_{k=1,\cdots,n}$.
\label{al:rhc}
\end{algorithm}
By {\it uniformly random} we mean that each components
of the initial models are chosen randomly with uniformly probability
between the maximum and minimum possible values. In this paper, all
RHC numerical results use non-linear Conjugate Gradient as transition
operators \cite{coool}.



\subsection{Landscape of Objective Functions}

Chavent \cite{curvature2} developed sufficient conditions for an
objective function to be locally convex. These conditions are based on
the distance $\times$ curvature induced by the objective function on
trajectories. In principle, this local convexity criterion could be
generalized to global samples of an objective function, to provide a
global measure of complexity.

On the other hand, imagine the surface of an objective function being
a high-dimensional landscape with hills and basins of different depths
and widths scattered on the surface. Performance of searching
algorithms depends to a large extent on topographical features on this
landscape.

For SA and GAs, this situation is summarized heuristically by Kaufmann
\cite{kaufmann2}: 
\begin{quote}
``Annealing works well only in landscapes in which deep energy wells
also drain wide basins. It does not work well on either a random
landscape or a ``golf course'' potential, which is flat everywhere
except for a unique ``hole''. In the latter case, the landscape offers
no clue to guide search.
\end{quote}
\begin{quote}
Recombination (in GAs) is useless on uncorrelated landscapes but
useful under two conditions (1) when the high peaks are near one
another and hence carry mutual information about their joint locations
in genotype space and (2) when parts of the evolving system are
quasi-independent of one another and hence can be interchanged with
modest chances that the recombined system had the advantage of both
parents.''
\end{quote}
In addition, since RHC uses local-descent transition operators, its
performance will also be strongly influenced by topography.

It has been proposed that functions can be characterized by their
spatial correlation properties \cite{weinberger,stadler2}. 
Several typical combinatorial optimization problems have been
investigated by studying the correlation in landscapes: the TSP
\cite{stadler4}, graph-bipartitioning problem \cite{stadler1}, and the
$NK$ model problems, a spin-glass like problem in biology 
\cite{kaufmannweinberger}. Using correlation features of the objective
function's landscape as a criterion, these authors study the
effectiveness of particular global algorithms for certain types of
landscapes. 

In addition, analyzing the topography of high-dimensional energy
functions is important in physics and chemistry. Berry and
Breitengraser-Kunz \cite{prl_topo} studied topography and dynamics of
multidimensional inter-atomic potential surfaces by analyzing a
population of local minima, each of which has two saddle points
connected to it. By connecting these samples in a certain order, the 
high-dimensional function surface is represented by a series of
one-dimensional lines. By looking at these one-dimensional plots, the
topography information is represented by the width and depth of the
primary, secondary or tertiary basins of attractions
\cite{prl_topo}.

The structure of high-dimensional Hamiltonians has also been 
studied by means of entropy \cite{prl_entropy}. For an
$N$-dimensional  Hamiltonian, a collection of local extrema are 
first found by
some means.  Contributions of these local-minima are
represented by a probability distribution $\{p_i\}_{i=1,\cdots,n}$
where 
\begin{equation}
 p_i \propto \Delta^{(i)}(N) = \prod_{k=1}^N \delta^i_k.
\label{eq:area}
\end{equation}
Here, $\delta^i_k$ is the estimated width of the $i$th basin of
attraction along the $k$th coordinate. The 
$N$-dimensional surface is then characterized by the following
entropy, 
\begin{equation}
 S(N) = - \sum_{i} p_i \ln p_i = \left \langle \ln
 \left (\frac{1}{\Delta(N)} \right ) \right \rangle.
\label{eq:entropy}
\end{equation}

In this paper, we use a similar measure.
However, we estimate $p_i$ by random
hill-climbing rather than equation~\ref{eq:area}. Further, we base our
measure not on $p_i$ itself, but rather on a related
probability that takes into account of the values of the local
minima. We also perform a confidence interval analysis.
Finally, as a concrete application, 
we show that this measure can be used to compute the optimal
simplification of a multi-resolution analysis (MRA) of highly non-convex
seismic optimization problem.

\section{Measures of Topography}

\subsection{Definitions}

The surface topography of functions is largely associated with
the number of local minima, widths of the basin of attractions
associated with these minima, and relative depths of these basins.
The {\it basin of
attraction} associated with the $i$th local minimum may be loosely
defined as {\it the maximum volume $A_i$ in the $N$-dimensional model
space within which all models can converge to the $i$th local minima
after infinite number of iterations by a local descent search
algorithm.} Suppose the volume of the entire model space is
represented as $M$, then the ratio  $p_i=\frac{A_i}{M}$ is the
probability of converging to the $i$th local minimum for a uniformly
random model.  The following definition serves to introduce
three quantitative measures of topography:  a probability associated
with the relative volumes of the basins of attraction ($p$), a
version of $p$ scaled by the estimated depths of the basins of
attractions ($q$) and the entropy of $q$.

\begin{definition} {\bf Entropy-Based Topography}\\
Let 
 $F: {\cal M} \subset {\bf R}^N \rightarrow{\cal Y}\subset {\bf R}$ be
bounded and have $n$ isolated local minima, where $n$ is finite. Let
 $\{{\bf m}_i\}_{i=1,\cdots,n}$ be these distinct local minima, and
 $\{y_i = F({\bf  m}_i)\}_{i=1,\cdots,n}$ be their corresponding
function values. Let $p_i$ be the probability that a model chosen with
uniform probability in ${\cal M}$ will converge to the $i$th local
minimum under the action of an exact local optimization algorithm.

Define $\{q_i\}_{i=1,\dots,n}$ to be a probability distribution where
\begin{eqnarray}
q_i \propto 
\begin{cases} 
p_i, & \hbox{if} \; \sigma = 0; \cr
p_i \, e^{-\frac{|y_i-y_m|}{\sigma}}, &\hbox{otherwise}, \cr
\end{cases}
\label{eq:prob1}
\end{eqnarray}
where $i \in [1,n]$, $y_m$ is the value at the global minimum,
 $\sigma = \frac{1}{n} \sum_{j=1}^n |y_j-y_m|$, and 
 $\sum_{i=1}^n q_i = 1$. 
The entropy is defined to be,
\begin{eqnarray}
   C_e & = & -\sum_{i=1}^n q_i \ln(q_i)
   \nonumber \\
        & = & \left\langle \ln \left(\frac{1}{q_i}\right)
        \right\rangle,
 \label{eq:measure}
 \end{eqnarray}
where the angle brackets denote the expected value with respect to the
probability distribution $\{q_i\}_{i=1,\dots,n}$.
\label{def:comp}
\end{definition}
For the entropy in Definition~\ref{def:comp}, it is always
the case that $C_e \ge 0$. If the function is unimodal (only one
extremum), then $C_e = 0$. On the other hand, if the function has $n$
isolated and equally valued local extrema, then $C_e$ is 
 $\ln n$. Therefore, the entropy increases with the number of local
extrema $n$.

\begin{figure}
\centerline{\includegraphics[width=100mm]{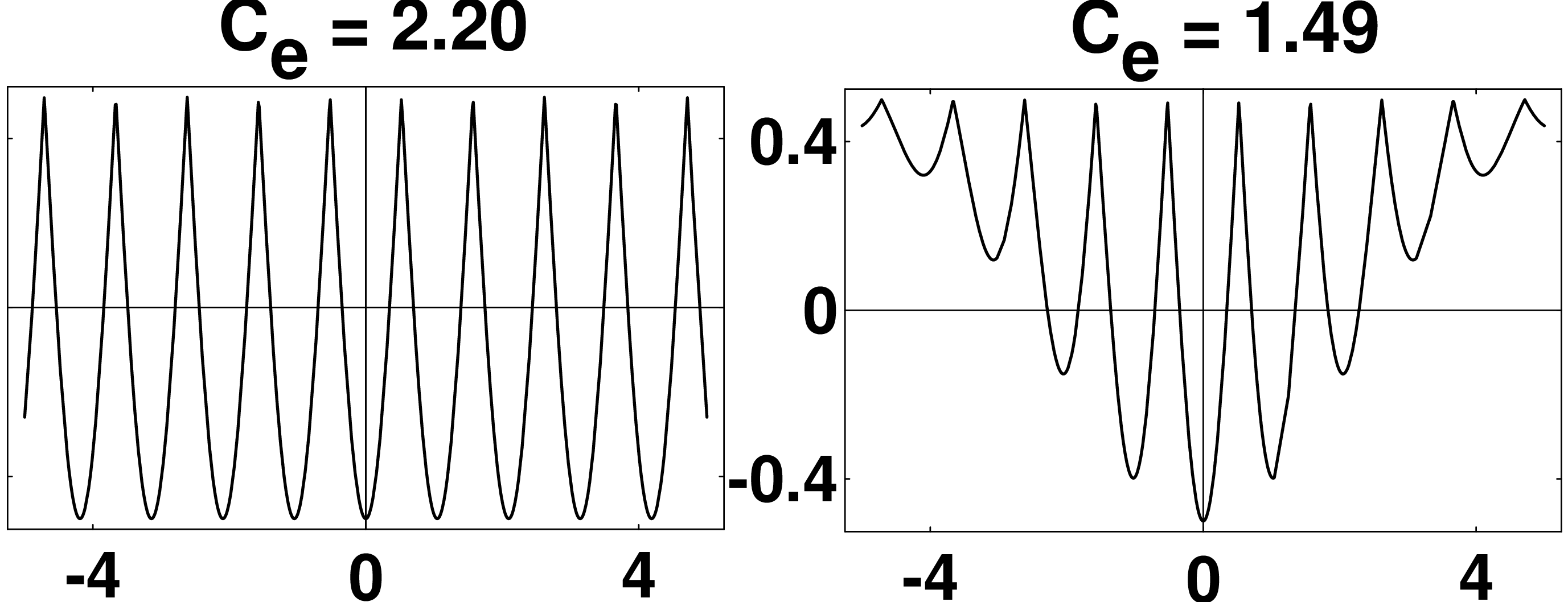}}
\caption{Two one-dimensional functions with the same number and
widths of basins of attractions. The function on the left have
entropy $C_e=2.2$, while the function on the right has
$C_e=1.5$.}
\label{f:cosines}
\end{figure}

As a simple example, Figure~\ref{f:cosines} shows two one-dimensional
functions with the same number of local minima and widths of basins of
attractions ($A_i \propto p_i$). However, the difficulty of
minimizing these functions is different: the left function has
identical basins of attractions, while the one on the right has a
dominant global minimum at $x=0$ and decreasingly important local
minima away from the center. The entropy of
Definition~\ref{def:comp} gives a higher $C_e$ value to that of the
function on the left ($C_e=2.2$) than that on the right ($C_e=1.5$).

Since the functions we are interested in can usually be evaluated only
point-wise, the number of local minima $n$ and 
$\{p_i\}_{i=1.\dots,n}$ are not known. Some degree of global sampling
is essential in order to achieve the characterization we seek. As
shown in Algorithm~\ref{al:rhc}, RHC explores various regions of the
model space and takes initial samples down-hill to the bottom of the
basins on the surface of functions. Therefore, a statistical analysis
of results of systematic RHC searches can be used to estimate the
topographic quantities.

Suppose the local-descent search is ideal, i.e. all initial models
converge to exact local minima, the number of models converging to
each local minima from $K$ randomly chosen initial models has a
multinomial probability distribution. If the initial models are
randomly chosen under a uniform probability distribution, the
probability of converging to the $i$th local minimum is proportional
to the width of the $i$th basin of attraction, $p_i$. 
Let $K_1, \dots, K_n$ be the random variables representing the
frequency of models converging to each of the local minima. For the
population of $K$, the joint probability density of these random
variables is 
\begin{equation}
f(k_1,\dots,k_n) = \frac{K!}{k_1!,\dots,k_n!}\; p_1^{k_1} \dots
  p_n^{k_n},
\label{eq:multinomial}
\end{equation}
where $\sum_{i=1}^{n} k_i = K$. For each $i \in [1, n]$, the mean
value of the random variable $k_i$ is $E[k_i] = K \, p_i$. 
Therefore, the convergence frequency of a RHC of large population can
be used to estimate the number of local minima as well as the widths
of basin of attractions. Hence, the estimation of the entropy
measure in Definition~\ref{def:comp} can be defined as follows.

\begin{definition} {\bf Entropy-Based Estimates} \\
Let $F:{\cal M}\subset {\bf R}^N \rightarrow{\cal Y}\subset {\bf R}$
be bounded and have finite number of isolated local minima.
Let $\{{\bf m}_i\}_{i=1,\cdots,\hat{n}}={\bf RHC}(F,K,\epsilon,cmax)$
be the distinct converged models of RHC searches.
Let $\{k_i\}_{i=1,\cdots,\hat{n}}$ be the frequency distribution of
the final population, and  
 $\{y_i = F({\bf  m}_i)\}_{i=1,\cdots,\hat{n}}$ be their corresponding
function values.

Define the estimated entropy ${\hat C}_e$ as
\begin{eqnarray}
{\hat C}_e & = & -\sum_{i=1}^{\hat n} {\hat q}_i \ln({\hat q}_i)
\nonumber \\ 
 & = & \left\langle \ln \left(\frac{1}{{\hat q}_i} \right)
        \right\rangle,
\label{eq:estimates}
\end{eqnarray}

where  ${\hat q}_i$ is normalized to a probability distribution
 \( \sum_{i=1}^{\hat n} {\hat q}_i = 1,  \) and
\({\hat q}_i \propto x_i\,v_i,\)
in which
\begin{eqnarray}
  x_i & \equiv & \frac{k_i}{K},
\label{eq:x}
\end{eqnarray}
and
\begin{eqnarray}
   v_i \equiv 
\begin{cases} 1, &\hbox{if} \; \sigma = 0; \cr
     \exp(-\frac{|y_i-y_{m}|}{\sigma}),
     &\hbox{otherwise}, \cr
\end{cases}
\label{eq:boltz}
\end{eqnarray}
where $y_m= \min{\{y_i\}}$ and
 $\sigma = \frac{1}{{\hat n}} \sum_{j=1}^{\hat n} |y_j-y_m|$.
\label{def:estimates}
\end{definition}

Definition~\ref{def:estimates} is a statistical estimation of the
entropy in Definition~\ref{def:comp}. The exact
entropy $C_e$ characterizes topographical features of objective 
function, and hence independent of numerical computation and any
searching technique. The estimation ${\hat C}_e$, however, would be
influenced by numerical issues.
If, for examples, the curvature of the function is nearly zero, which
is equivalent to an ill-conditioned Hessian matrix, gradient-based
local descent searches may not converge to the exact local minima. The
estimated value of ${\hat C}_e$ in such a situation may be
higher than the true complexity $C_e$.
In practice, however, it is often difficult to distinguish the results
of such ill-conditioning from those of multi-modality.
Therefore, taking such numerical issues into account can represent
an important aspect in the difficulty of optimization. 

\subsection{Numerical Examples}

In this section, we use the entropy ${\hat C}_e$ to study two
commonly used test functions in optimization, the Rosenbrock and
Griewank functions.

\subsubsection{$N$-dimensional Rosenbrock function}

\begin{figure}
\centerline{\includegraphics[width=42mm]{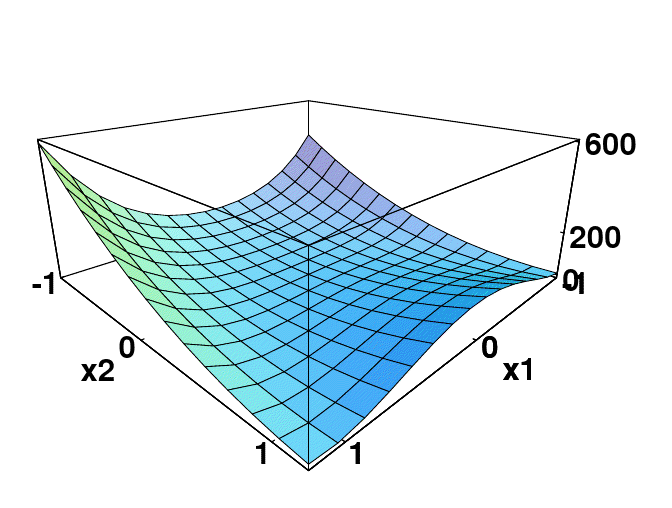}
\includegraphics[width=42mm]{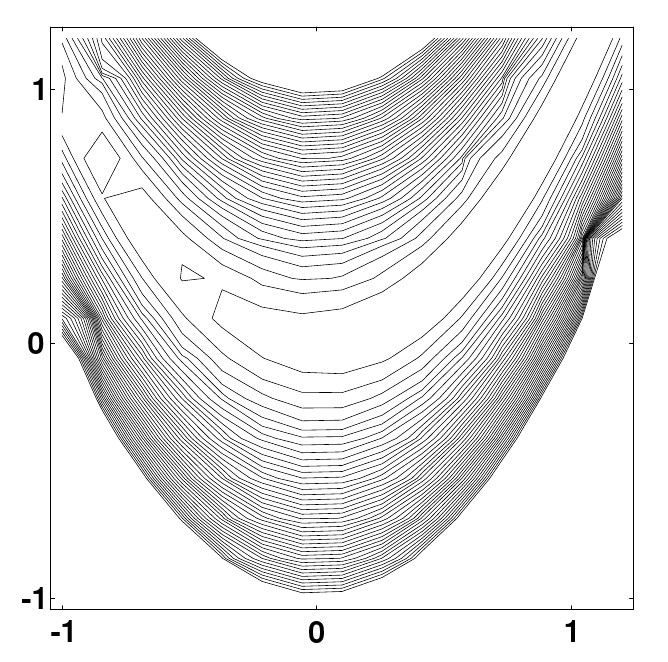}}
\caption{Two-dimensional Rosenbrock function. The figure on the left
  is a 3-D plot of the function surface, while the one on the right
  shows the contour plot of the same function.}
\label{f:rosen2d}
\end{figure}

An $N$-dimensional Rosenbrock function can be written as
\begin{equation}
R({\bf x}) = \sum_{i=1}^{N-1} \left [ 100 (x_i - x_{i-1}^2)^2 +
  (1-x_{i-1})^2 \right ],
\label{eq:rosenbrock}
\end{equation}
where ${\bf x} = (x_0, \; \dots, \; x_N)$.
Although unimodal, the long and narrow basin is a challenge
for searching algorithms.  Figure~\ref{f:rosen2d} shows the
function surface and its contour when $N=2$. When $N \ge 2$, the
function is still unimodal, but it is not easy to see how the
increase of dimensionality alters the difficulty of optimization.

\begin{figure}
\centerline{\includegraphics[width=85mm]{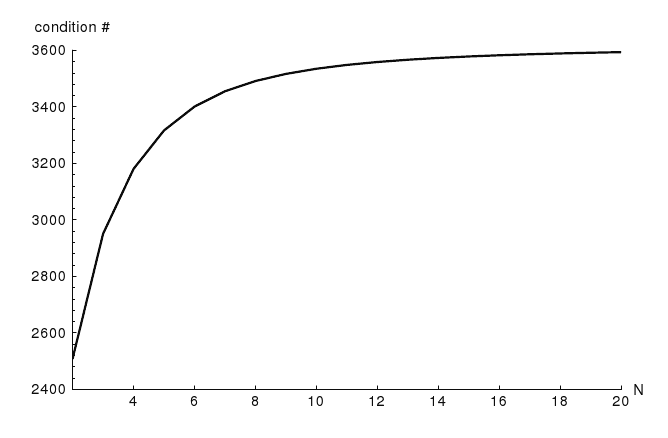}}
\caption{Condition number of the Hessian matrix for the N-dimensional
  Rosenbrock function at the global minimum.}
\label{f:conditions}
\end{figure}

\begin{figure}
\centerline{\includegraphics[width=85mm]{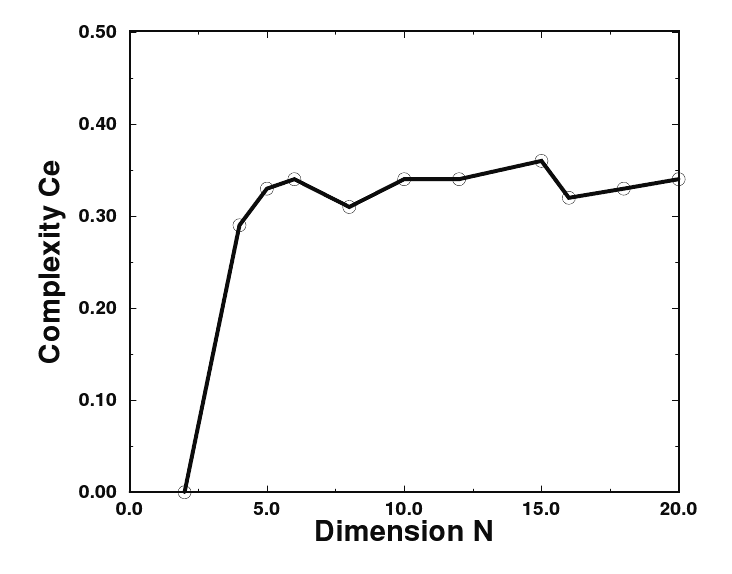}}
\caption{${\hat C}_e$ for N-dimensional Rosenbrock
functions as a function of $N$.}
\label{f:comp_rosen}
\end{figure}

One way of studying the spatial curvature of functions is by looking
at the ratio of largest and smallest eigenvalues ({\it condition number})
of the Hessian at a point. The Hessian for
equation~(\ref{eq:rosenbrock}) is a tri-diagonal matrix,
\begin{eqnarray}
\left ( \begin{array}{ccccc}
 a_0 & c_0 & 0 & \cdots & 0 \\
b_1 & a_1 & c_1 & 0 & \cdots \\
0 & b_2 & a_2 & c_2 & 0 \\
\vdots & \vdots & \vdots & \vdots & \vdots \\
 0 & 0 & \cdots & b_{N-1} &  a_{N-1}
\end{array}  \right)
\label{eq:hessian}
\end{eqnarray}
where 
\begin{eqnarray*}
a_i & = & 
\begin{cases} 
2 + 1200\,x_0^2-400\,x_1, & \hbox{if} \; i = 0; \cr
 202 + 1200x_{i-1}^2 - 400x_i, & \hbox{if}  $0 < i < N-1\nonumber$; \cr
 200, & \hbox{if} i = N-1, \cr
\end{cases}
\end{eqnarray*}
\begin{eqnarray*}
 b_i & = & -400\,x_{i-1}, \;\;\;\;  0 < i \le N-1, \\
 c_i & = & -400\,x_{i+1}, \;\;\;\;  0 \le i < N-1.
\end{eqnarray*}
At the global minimum $(1,\,1,\,\cdots,\,1)$, the tri-diagonal matrix
equation~(\ref{eq:hessian}) becomes Toeplitz except for $a_0$ and
 $a_{N-1}$.
The condition number of the Hessian at the global minimum reaches an
asymptote with increasing dimension, as shown in
Figure~\ref{f:conditions}.
Figure~\ref{f:comp_rosen} shows ${\hat C}_e$ as
a function of the number of dimensions; it shows the same asymptotic
trend as does the condition number. Thus the increasing complexity for
low dimensions is the result of increasing ill-conditioning of the
Hessian and has nothing to do with local minima.

\subsubsection{High-dimensional Griewank functions}

\begin{figure}
\centerline{\includegraphics[width=100mm]{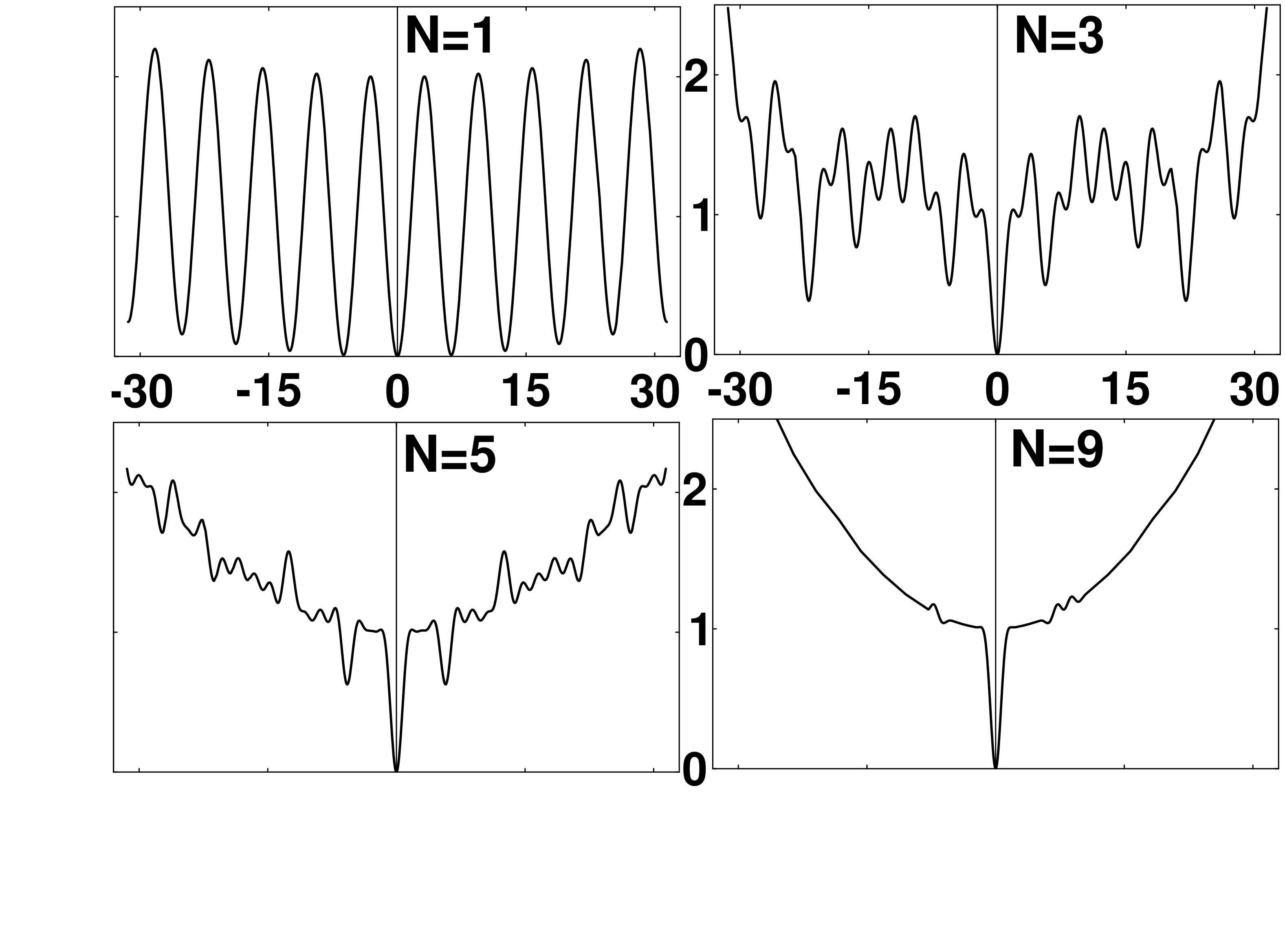}}
\caption{Diagonal slices of $N$-dimensional Griewank functions.}
\label{f:griew-diag}
\end{figure}

The Griewank function is also used to test optimization algorithms
 \cite{whitley_testfcn2,torn-zilinska}:
\begin{equation}
g({\bf x}) = 1 + \sum_{i=1}^N \frac{x_i^2}{4000} - \prod_{i=1}^N
\cos \left(\frac{x_i}{\sqrt{i}} \right)
\label{eq:griew}
\end{equation}
The cosine term makes equation~(\ref{eq:griew}) multi-modal.
Figure~\ref{f:griew-diag} shows a one-dimensional slice of the
Griewank function along the diagonal of the hypercube for dimensions 
 $1,3,5,9$. Whitley et al. \cite{whitley_testfcn2} observed such
 slices and concluded that ``as the dimensionality increases the 
local optima induced by the cosine decrease in number and
complexity''.

\begin{figure}
\centerline{\includegraphics[width=100mm]{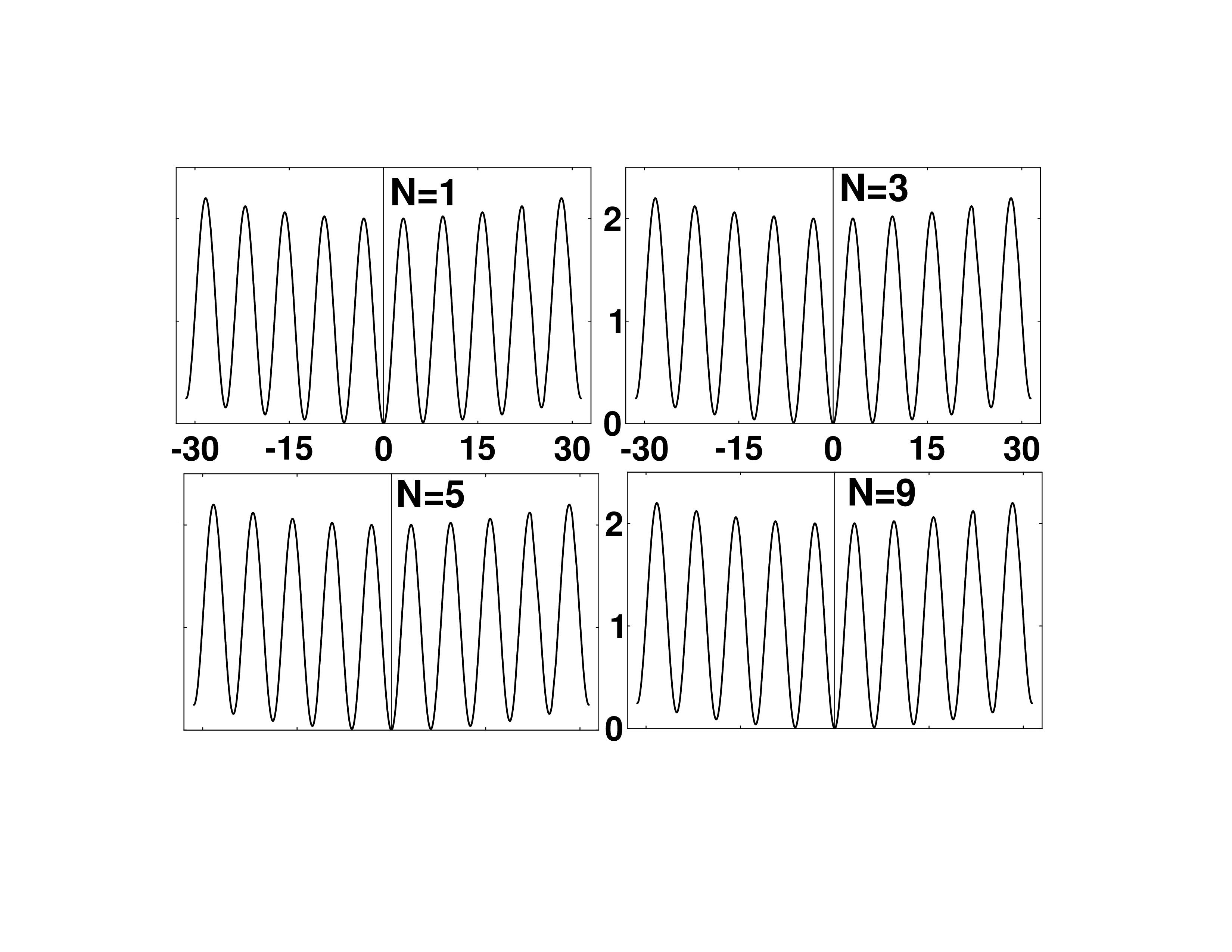}}
\caption{Slices of $N$-dimensional Griewank functions. All variables
  but one are fixed at $0$.}
\label{f:griew-orig}
\end{figure}

\begin{figure}
\centerline{\includegraphics[width=100mm]{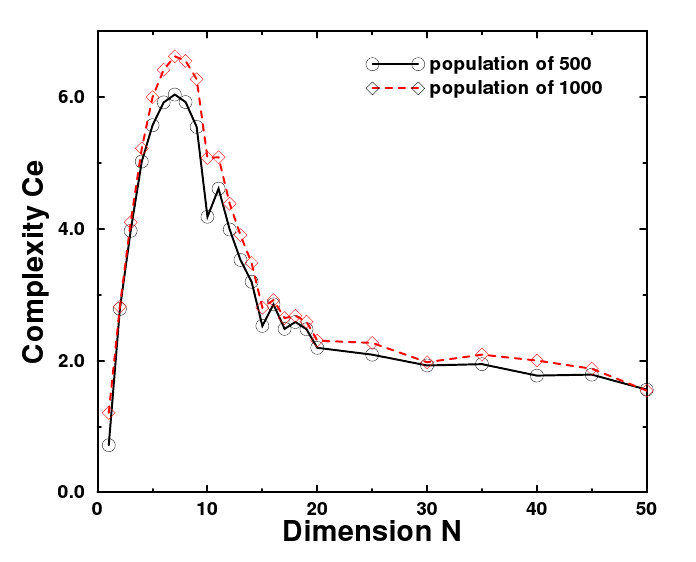}}
\caption{$\hat {C_e}$ as a function of dimension $N$ for the 
  Griewank function with populations of $500$ and $1000$.}
\label{f:griew}
\end{figure}

However, such pictures can be misleading since they tell us only about 
low-dimensional projections of the function. Figure~\ref{f:griew-orig}
shows slices of the same functions when all but one variables are
fixed to be $0$. The increasing dimensionality does not change the
oscillation around the global minimum at the origin.
Therefore, studying the overall performance of high dimensional
functions could be tricky. We compute ${\hat C}_e$ for the Griewank
function with a population $500$ and $1000$ models in the hyper-cube
of  $-10 \le x_i \le 10, \;\;\; i=0,\cdots,N-1$. Figure~\ref{f:griew}
shows the resulting ${\hat C}_e$ for dimensions up to $50$
for initial populations of both $500$ and $1000$. Both curves in
Figure~\ref{f:griew} give us consistent results that the complexity of
Griewank function in this range increases till dimension around $9$,
then decreases when number of dimension continuous to increase. This
result can be verified by the analysis of Griewank
function. Therefore, using the entropy we can
understand more comprehensively
the dimensional-dependence of complexity
of certain functions than by simply looking at hyper-planes.

\section{Confidence Intervals Analysis}

Next we derive confidence intervals on the entropy in
Definition~\ref{def:estimates}. The following analysis is based on the
assumption of {\it ideal RHC}, which is a special case of
Algorithm~\ref{al:rhc} where an infinitely large $cmax$ is allowed and
 $\epsilon$ is infinitely small.
First, it is easy to prove that as long as the population size $K$ is
large enough, $x_i$ defined in equation~(\ref{eq:x}) would be good
approximation to $p_i$ for $\forall i \in [1,n]$. 
We have the following theorem, the proof of which is given in the
appendix.

\begin{theorem}
Let 
 $F: {\cal M} \subset {\bf R}^N\rightarrow {\cal Y}\subset {\bf R}$
be bounded and have finite number of isolated local minima. 
Let $p_i$ be the probability of converging to the $i$th local minimum
for a starting model chosen with uniform probability on ${\cal M}$.
Perform an ideal RHC as defined in Algorithm~\ref{al:rhc} with an
initial population of $K$.
Let $\alpha$ and $\beta$ be related by the following equation, 
\begin{equation}
Pr(|z| \le \beta) = 1 - \alpha,
\label{eq:alpha-beta}
\end{equation}
where $z$ has a standard-normal distribution $N(0,1)$.

Let $x_i$ be as defined in equation~(\ref{eq:x}).
If the population $K$ is such that $K\, p_i \ge 5$ for 
any $i \in [1,n]$, we have the following,

\begin{enumerate}
\item $x_i$ has an approximate normal distribution with 
\begin{equation}
E[x_i] = p_i,
\label{eq:x-mean}
\end{equation}
and
\begin{equation}
var(x_i) = \frac{p_i\, (1-p_i)}{K}.
\label{eq:x-variance}
\end{equation}
\item $x_i$ is an unbiased, consistent estimator of $p_i$.
\item With confidence of $(1-\alpha) \%$, the error associated with
estimating $p_i$ from $x_i$ is bounded by 
\begin{equation}
 \beta \, \sqrt{\frac{x_i (1-x_i)}{K}}.
\end{equation}
\end{enumerate}
\label{th:intervx}
\end{theorem}

To get some ideas of the magnitudes of the population size and the
confidence interval, here is a simple example.
\begin{example}
If for a problem as described in Definition~\ref{def:estimates}, we
have $p_i = 0.01$. Then for approximating the binomial distribution
with a normal distribution, we need at least $K > 500$.
\label{ex:first}
\end{example}
\begin{example}
For the same problem as stated in Example~\ref{ex:first}, suppose a
population size of $520$ was used in an ideal RHC, and $100$ of the
models converged to the $i$th local minimum. Then, 
 $x_i \approx 0.192$. If we want to have $90\%$ confidence, then
 $\beta = 1.65$. The error bound for the estimation of $p_i$ with
 $x_i$ would be $0.0285$. That is, with $90 \%$ confidence, we can say
that  $0.175 \le p_i \le 0.220$.
If, on the other hand, we want to have $95 \%$ confidence for this
estimation when $\beta = 1.96$, then $0.158 \le p_i \le 0.226$.
\end{example}

In the following theorem, we estimate the distribution and error bound
of the estimation for the entropy
Definition~\ref{def:estimates}. The proof of the following theorem is
also given in the appendix.
\begin{theorem}
Let 
 $F: {\cal M} \subset {\bf R}^N\rightarrow{\cal Y} \subset {\bf R}$ be
bounded and have finite number of isolated local minima.
Let $p=\{p_i\}_{i=1,\dots,n}$ be the probability distribution in
Definition~\ref{def:comp}, and $p_m=\min\{p_i\}$.
Let $C_e$ be the entropy of $F$ (as in Definition~\ref{def:comp})
and suppose the RHC of population $K$ is ideal. As a result, the
initial population of models converge to different local extrema with
a frequency distribution of $\{k_i\}_{i=1,\dots,{\hat n}}$. Finally,
let $x_i$ be defined as in equation~\ref{eq:x}, and 
 ${\hat C}(x_1, x_2,\dots, x_{\hat n})$ the estimated entropy (as
in Definition~\ref{def:estimates}). If $\alpha$ and  $\beta$ are
defined as in equation~(\ref{eq:alpha-beta}), then we can prove the
following statements:

\begin{enumerate}
\item 
${\hat C}_e$ has an approximate normal
distribution with 
\begin{equation}
E[{\hat C}_e] = C_e,
\label{eq:ce_mean}
\end{equation}
and
\begin{equation}
var({\hat C}_e) = 
  \sum_{i=1}^n \frac{c\,v_i \, (1+\ln q_i)^2}{K}\,
 q_i - \frac{(1- C_e)^2}{K},
\label{eq:ce_variance}
\end{equation}
in which
\begin{equation}
 q_i = c \, v_i \, p_i,
\end{equation}
for $\forall i \in [1, n]$, $c$ is a scale factor so that
 $\sum_{i=1}^n q_i =1$, and  $v_i$ is as defined in
equation~(\ref{eq:boltz}). 

\item
 ${\hat C}_e$ is an unbiased, consistent estimator of $C_e$.

\item 
Let ${\hat q}_{min} = \min_{i\in [1,{\hat n}]}\{{\hat q}_i\}$.
If the population size $K$ is such that $K\,p_m \ge 5$, then with
confidence of $(1-\alpha)\,100\%$, the estimation error of the
complexity is at most  $\beta\,\delta$ where $\delta \ge 0$. That is,
\begin{eqnarray*}
| C_e - {\hat C}_e (x_1,\dots, x_{\hat n})| \le \beta\,\delta,
\end{eqnarray*}
where
\begin{equation}
  \delta^2 = \frac{c}{K} (1 - (2 + \ln {\hat q}_{min}) \,
  {\hat C}_e)  - \frac{(1 - {\hat C}_e)^2}{K} > 0.
\label{eq:delta}
\end{equation}

\end{enumerate}
\label{th:confidence}
\end{theorem}

\begin{remark}
When the number of local minima, $n$, is large, the real $p_m$ may
very small. Then, an unrealistically large initial population size $K$
may be required to satisfy $K \, p_m \, \ge \, 5$. Realistically, we
have to content ourselves with not being able to find all local minima
in such difficult situations. If the smallest basin we found with a
$K$-population RHC is  ${\hat p}_m$ and ${\hat p}_m > p_m$, there are
some $p_i < \{\hat p\}_m$ which are not found by the RHC. Their
corresponding $q_i$ would not be accounted for in the complexity
estimation. 
However, the contributions of these narrow basins to the complexity
are proportional to $q_i \log q_i$. Since, 
 $\lim_{q  \rightarrow 0} q_i \log q_i = 0$, the error caused by these
narrow basins will be small as long as $q_i = c\, v_i \, p_i$ is
small. Since $0 < v_i \le 1$ for $\forall i \in [0,n]$ by definition,
these conditions can be easily satisfied as long as these narrow
basins are not global minima.
\end{remark}

We conclude this section by showing an example of the evaluation of
the confidence interval for the Griewank function. 
It is important to note that in such an analysis, it is assumed that
the RHC algorithms are exact. That is, numerical effects are ignored.

\begin{example}
We want to evaluate the confidence interval for the complexity
calculation of  the $9$-dimensional Griewank function shown in
equation~(\ref{eq:griew}). In Figure~\ref{f:griew}, we show that the
complexity estimation for the population of $K=1000$ is $6.27$.
Using the calculated data, $\ln q_{min}=-13.62$, $c=1.86$, and 
 ${\hat  C}_e=6.37$, we can estimate that $\delta \approx 0.33$. 
So, with $90 \%$ confidence, the error bound would be 
 $\pm 0.65$. Therefore, we can say that the true complexity value
 $C_e$ is between $6.9$ and $5.6$.
\end{example}

\section{A Geophysical Application}

\subsection{Estimating Near-Surface Heterogeneities}

\begin{figure}
\centerline{\includegraphics[width=80mm]{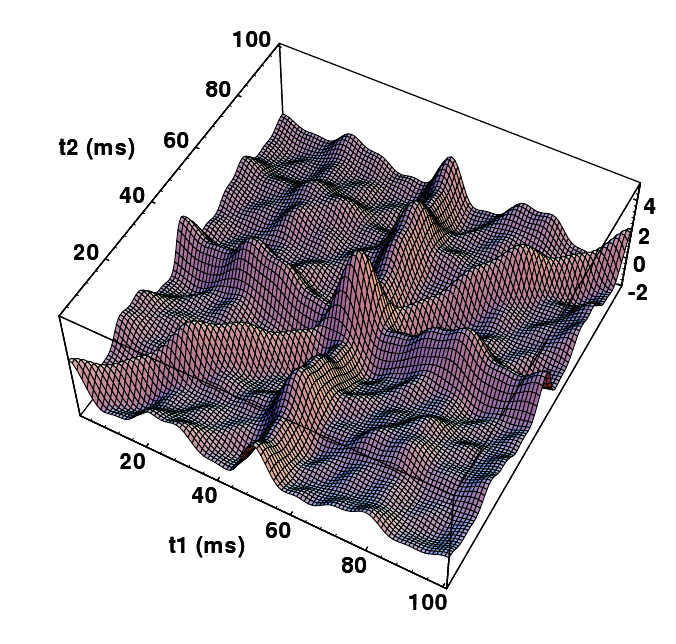}}
\caption{Landscape of a 2-D residual statics objective function.}
\label{f:data3D}
\end{figure}

In exploration seismology, ``statics'' are the time shifts in seismic
reflection data caused by heterogeneous material properties in the
near surface.  This causes jitter in the data and degrades processing
procedures designed to enhance signal-to-noise, such as averaging.
It is possible to formulate an optimization procedure for these
static time shifts (the objective function being the power of the
averaged data as a function of time shifts), 
but the resulting optimization problem is
highly non-convex \cite{rothmana}.  This is illustrated in
Figure~\ref{f:data3D} with a toy example.


Consider an example where we need to align three otherwise identical
traces.  Fixing the first trace, we look for time-shifts for the
second and third traces, $t_1, t_2$, so that the sum of squares of the
stacked traces (stacking-power) is maximized. Figure~\ref{f:data3D}
shows an example of such a two-dimensional objective function, which
has hills and basins of attractions scattered on the landscape. In
practice, however, the stacking-power objective function is
high-dimensional and highly multi-modal. Monte Carlo global
optimization have become an important tool for solving large-scale
statics problems \cite{rothman85,rothman86}.

\begin{figure}
\centerline{\includegraphics[width=80mm]{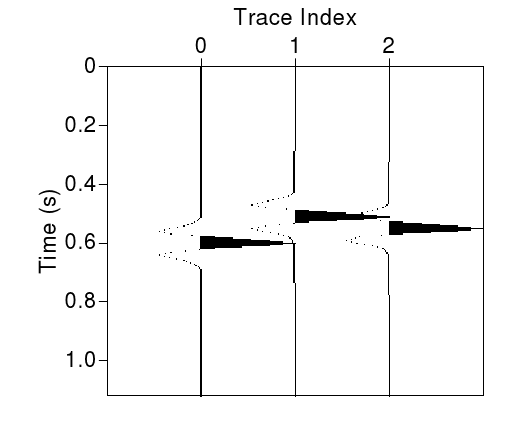}}
\caption{Three synthetic seismic traces that are shifted by
 random statics.}
\label{f:three}
\end{figure}

\begin{figure}
\centerline{\includegraphics[width=100mm]{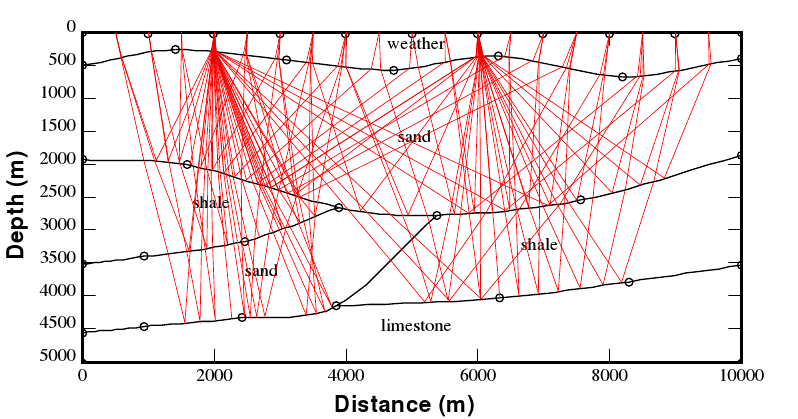}}
\caption{A hypothetical model of the Earth's upper crust showing rays
 associated with seismic waves propagating down from sources on the
 surface, reflecting off geologic boundaries and traveling upwards
 to receivers on the surface. Static-shifts of the seismic traces are
 caused by the combined time-distortions of near-source and
 near-receiver heterogeneities in the weathering layer.} 
\label{f:demo_stats}
\end{figure}

Figure~\ref{f:data3D} shows a statics objective function with two
unknowns. In practice, however, the time-shifts of the traces are not
independent. The statics of each trace are caused by the combined time
distortion of near-source and near-receiver heterogeneities ({\it
source-statics} and {\it receiver-statics}).
Figure~\ref{f:demo_stats} illustrates the similarity of travel paths
near each source and each receiver.

The recorded reflection seismic signals are usually sorted into
{\it midpoints} $y$ (of the source and receiver locations) and {\it
offsets} $h$ (half distance between the source and receivers).  Letting
$\bf{s}$ and $\bf{r}$ be unknown vectors of source- and
receiver-statics, this optimization problem can be formulated as
\begin{equation}
 \max_{\bf{s}, \bf{r}} F(\bf{s}, \bf{r}) = \sum_y{\sum_{h_1 \neq
     h_2} {\Phi^y_{h_1,h_2}(\tau(\bf{s}, \bf{r}))}},
\label{eq:cross}
\end{equation}
where $\Phi^y_{h_1,h_2}(\tau)$ is the cross-correlation between traces
(after a correction for propagation
effects known as ``normal move-out'' has been 
applied) of offsets $h_1$ and $h_2$ at midpoint $y$ evaluated at 
\[ \tau=s_{i(y,h_1)}+r_{j(y,h_1)}-s_{i(y,h_2)}-r_{j(y,h_2)}, \]
and $i(y,h)$ and $j(y,h)$ are the source and receiver indices for
midpoint $y$ and offset $h$, respectively. The function 
$F(\bf{s},\bf{r})$ in equation~(\ref{eq:cross}) is 
called the {\it stacking-power function}.

\begin{figure}
\centerline{\includegraphics[width=100mm]{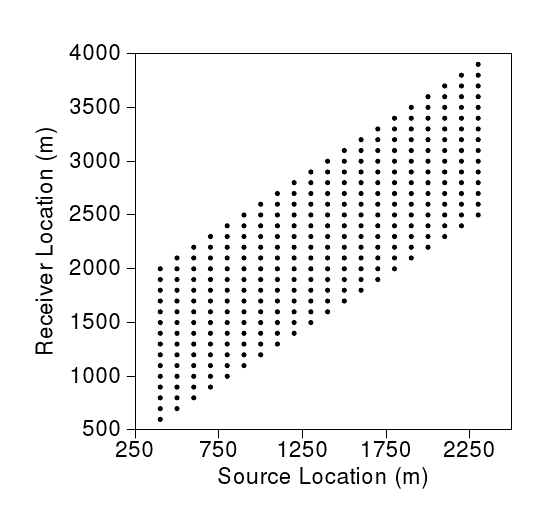}}
\caption{ 
Recording geometry of a synthetic data set. The horizontal axis is
the source position and vertical axis the receiver position.
}
\label{f:chart}
\end{figure}

\begin{figure}
\centerline{\includegraphics[width=100mm]{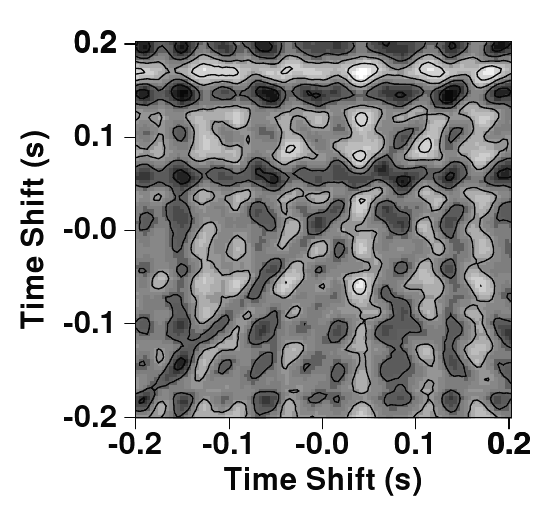}}
\caption{Hyperplane through the stacking-power
function for a problem with 55 unknown static time shifts.
}
\label{f:cont_no}
\end{figure}

Figure~\ref{f:chart} shows the recording geometry of one example
synthetic data set. This data set has $20$ sources, $35$ distinct
receivers and $320$ traces. All traces are identical except for random
source and receiver statics. These are generated by repeatedly
shifting  a single trace of field data. Thus, the objective
function of equation~(\ref{eq:cross}) has $55$ unknowns. When there
are no statics in the data, the global maximum of the function is at
the origin ($ {\bf  s}_i={\bf 0},\;{\bf  r}_j={\bf 0}$).
Figure~\ref{f:cont_no} shows an arbitrary 2-D hyper-planes of the
stacking-power function along the $10$th source and $20$ receiver
statics.

We have analyzed a realistic synthetic
statics problem involving some 320 seismic traces and 55
unknown static time shifts.  A hyperplane through the objective
function for this problem is shown in figure~\ref{f:cont_no}.
In addition to simply computing the entropy of this function
we will show how the entropy might be used to quantitatively
address issues related to the topography of functions.

\subsection{Behavior of the Multi-Resolution Analysis}

Rather than using Monte Carlo global optimization methods to solve
the statics problem as in \cite{rothman85,rothman86}, 
Deng \cite{lydiawavelet} has proposed, without proof, simplifying the
optimization via a multi-resolution analysis (MRA) of the seismic
traces. The idea is to use a wavelet decomposition to generate
successively simpler representations of the seismic data, thereby
eliminating progressively more local extrema from the objective
function. To be precise, let us define a Multi-Resolution RHC
algorithm: 
\begin{algorithm} {\bf MRHC}
 $({\bf {\vec m}} = {\mathit MRHC}(F,L,\epsilon,cmax))$ \\
Let $\{{\bf S}_i\}_{i=L,\cdots,0}$ be a sequence of decreasingly
smooth operators to be defined below, with ${\bf S}_0$ an identity
operator. 
\begin{enumerate}
\item Let $f_L = {\bf S}_L\, F$; choose an initial population
 $\{{\bf m}^0_k\}_{k=1,\cdots,K}$ with size $K$ at random;
apply Algorithm~\ref{al:rhc}, so
 $\{{\bf m}_k\}_{k=1,\cdots,M_L} =  {\mathit RHC}(f_L,K,\epsilon,cmax)$,
 and $i=L-1$.
\item Let $f_i ={\bf S}_i\, F, \; (L > i \ge 0)$ and 
  ${\bf {\vec m}}^0 = \{{\bf m}_k\}_{k=1,\cdots,M_{i-1}}$; run 
Algorithm~\ref{al:rhc},  
 $\{{\bf m}_k\}_{k=1,\cdots,M_i} = {\mathit  RHC}(f_i,K,\epsilon,cmax)$.
\item Decrease the level index $i$ by $1$, repeat $2$ until $i=0$. The
  final set of   models ${\bf {\vec m}}$ is the solution.
\end{enumerate}
\end{algorithm}

The smoothing operators $\{{\bf S}_i\}_{i=L,\cdots,0}$ could be a
sequence of low-pass filters with increasingly wider pass-band
\cite{bunks-saleck}, or a sequence of increasingly fine wavelet
operators \cite{lydiawavelet} for decomposing the input seismic data.
The sequence of smoothing operators should be such that the resulting
functions, $\{f_i\}_{i=L,\cdots,0}$, have the same global feature as
does the objective function $F$ for all levels and have decreasing
number of local optima when the level increases, and $f_0 = F$.
Deng \cite{lydiawavelet} showed that this could be achieved using the
shift-invariant wavelet basis of Saito and Beylkin \cite{SA-BEY}.

\begin{figure}
\centerline{\includegraphics[width=100mm]{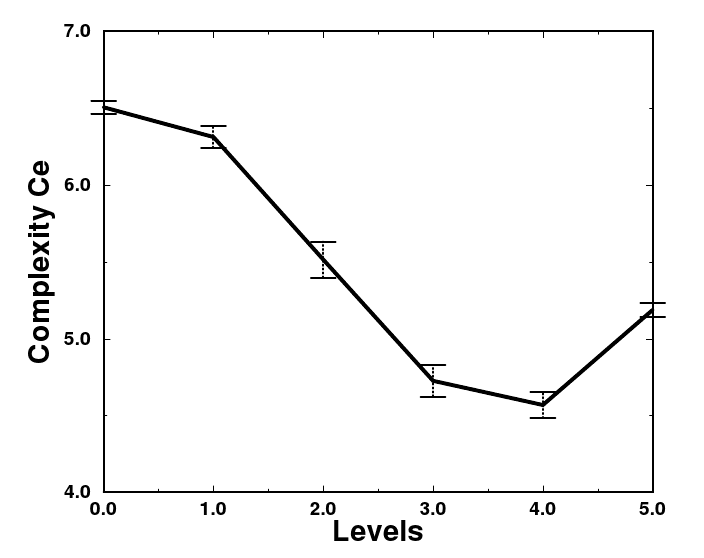}}
\caption{Complexity ${\hat C}_e$ as a function of the level of wavelet
  decomposition.}
\label{f:comp_stats}
\end{figure}

We now apply the entropy-based estimation of complexity to study
the multi-resolution analysis of the 55 parameter statics problem
introduced in the previous section.
Figure~\ref{f:comp_stats} shows ${\hat C}_e$ as a function of the wavelet
decomposition level using population $K=1000$; the mean and
one-standard deviation error bars are obtained from 32 independent
calculations.   Results are shown for 6 levels of decomposition
using a wavelet operator $\{ {\bf S}\}_{i=0,...,5}$ where $i=0$
is an identity operator, corresponding to use of the original data.
These results indicate that for this particular problem a complexity 
minimum is achieved for a wavelet decomposition of level 4.
Higher levels of decomposition actually increase the
complexity; presumably this results from the objective function
being too flat for local optimization.
Thus, the complexity measure gives us a way of choosing
a wavelet decomposition level to achieve optimal simplification
of an objective function.

\section{Conclusions}

We have developed a collection of simple
tools for analysis of the
topographic complexity of functions based on the
application of local optimization to randomly chosen
starting models.  In particular
we estimate
the number of basins of attractions on the function
landscape, the widths and depths of these basins
and the entropy of the resulting probabilities.

Assuming local descent searches are ideal, we have computed the
confidence intervals for the sampling error associated with this
complexity measure. There are, on the other hand, several practical
issues that we have neglected in this error analysis. Among them, we
can mention the convergence error 
caused by the finite computing time and the finite precision of the
local descent algorithms, the criterion for clustering of converged
models, and the size of the assumed smallest basin of attraction,
 $p_m$. These issues can be investigated by a Monte Carlo analysis as
shown in Figure~\ref{f:comp_stats}. 

\section{Acknowledgments}
This work is dedicated to the memory of Albert Tarantola.
The authors thank Dr. Bill Navidi for useful discussions and comments
on a draft of this work.  
This work was begun while the authors were at the Center for Wave Phenomena.

\appendix

\section{Proofs of the Confidence Interval Analysis}

\subsection{Proof of Theorem\ref{th:intervx}}

{\bf
Proof:} \\

Let $K_1, \dots, K_n$ be the random variables that represent the
frequency of initial models converging to each local minima.
The joint probability density for random variables
 $K_i, \; i \in [1,n]$ for population of $K$ is a multinomial
distribution. The marginal distribution for each of the random
variables is,
\begin{equation}
 f(k_i) = \frac{K!}{k_i ! (K-k_i)!} p_i^{k_i} (1-p_i)^{K-k_i},
\end{equation}
and the corresponding statistical quantities are,
\[
E[k_i]  =  K \, p_i, \;\;\;\; var(k_i)  =  K p_i (1-p_i),
\]
for $\forall i \in [1, n]$.

\begin{enumerate}

\item
For $K$ such that $K\, p_i \ge 5$, the above binomial distribution can
be approximated by a normal distribution. That is, the random
variable,
\begin{eqnarray*}
 Z_i & = & \frac{k_i - K p_i}{\sqrt{K p_i (1-p_i)}} \\
  & = & \frac{x_i - p_i}{\sqrt{p_i (1-p_i) / K}}
\end{eqnarray*}
approaches to standard normal distributions $N(0, 1)$
 (Theorem~{\bf 6.8} of \cite{freund}), where $x_i$ is
defined in equation~(\ref{eq:x}). Therefore, $x_i$ has a normal
distribution with the mean and variance as in
equations~(\ref{eq:x-mean}) and (\ref{eq:x-variance}).

\item 
From equation~(\ref{eq:x-mean}), we see that $x_i$ is an unbiased
estimator of $p_i$. Since $var(x_i) \propto \frac{1}{K}$,  we have
\begin{equation}
\lim_{K \rightarrow \infty} var(x_i) = 0.
\end{equation}
Therefore, 
\[  \lim_{K \rightarrow \infty} x_i = p_i. \]
 $x_i$ is also a consistent estimator of $p_i$ for
each $i \in [1,n]$.

\item
Now with confidence of $(1-\alpha) 100 \%$, we have 
\begin{equation}
  |p_i - x_i |  \le \beta \sqrt{\frac{p_i \, (1-p_i)}{K}},
\label{eq:interval}
\end{equation}
where the value of  $\beta = z_{\alpha/2}$ can be looked up from a
standard normal distribution table.

However, we do not know $p_i$ in advance. Approximating $p_i$
by $x_i$ when $K$ is large, we have the confidence interval for
the true $p_i$
\begin{eqnarray}
   |p_i - x_i|  \le \beta \sqrt{\frac{x_i \, (1-x_i)}{K}}
\label{eq:intervalx}
\end{eqnarray}
for $i \in [1, n]$ (Theorem~{\bf 11.6} of \cite{freund}).

\end{enumerate}

{\Large $\Box  $}

\subsection{Proof of Theorem~\ref{th:confidence}}

{\bf Proof:} \\

From Theorem~\ref{th:intervx}, we know that each random variable $x_i$
for $i \in [1, n]$ has an approximate normal distribution
 $N(p_i, \frac{p_i\, (1-p_i)}{K})$ when the population $K$ is such
that $K\, p_i \ge 5 $.
Since ${\hat q}_i = c v_i x_i$, then  ${\hat q}_i$ also has an
approximate normal distribution
 $N(c\,v_i\,p_i, \frac{c^2\, v^2_i \, p_i\, (1-p_i)}{K})$.
Since ${\hat q}_i$ would be very close to $q_i$ when $K$ is large, we
can make the following approximation,
\begin{equation}
{\hat q}_i \ln {\hat q}_i \approx q_i \ln q_i + (1 + \ln q_i)({\hat
  q}_i - q_i),
\label{eq:approx}
\end{equation}
which is a linear function of the random variable ${\hat q}_i$.
Therefore, ${\hat q}_i \ln {\hat q}_i$ is also approximate normal
distribution,
\begin{eqnarray}
 E[{\hat q}_i \ln {\hat q}_i] & = & q_i \ln q_i, \\
 var({\hat q}_i \ln {\hat q}_i) & = & (1 + \ln q_i)^2
\frac{c^2\, v^2_i \,p_i\, (1-p_i)}{K}.
\label{eq:mean}
\end{eqnarray}

\begin{enumerate}

\item
Since ${\hat C}_e$ is a linear combination of 
 ${\hat q}_i \ln {\hat  q}_i$ for $i \in [1, n]$, 
 ${\hat  C}_e$ also has an approximate normal distribution.
Then, we have,
\begin{eqnarray*}
E[{\hat C}_e] = - \sum_{i=1}^n E[{\hat q}_i \ln {\hat q}_i]
  = - \sum_{i=1}^n  q_i \ln q_i = C_e,
\end{eqnarray*}
and
\begin{eqnarray*}
var({\hat C}_e) = \sum_{i=1}^n var({\hat q}_i \ln {\hat q}_i) +
  \sum_{i \ne j} Cov({\hat q}_i \ln {\hat q}_i, {\hat q}_j 
 \ln {\hat q}_j).
\end{eqnarray*}

For calculating 
 $Cov({\hat q}_i \ln {\hat q}_i, {\hat q}_j \ln {\hat q}_j)$, recall
equations~(\ref{eq:approx}) and (\ref{eq:mean}), 
\begin{eqnarray*}
Cov({\hat q}_i \ln {\hat q}_i, {\hat q}_j \ln {\hat q}_j) & = &
 E[({\hat q}_i \ln {\hat q}_i - q_i \ln q_i) ({\hat q}_j \ln {\hat
   q}_j - q_j \ln q_j)] \\
 & = & (1+\ln q_i)\,(1+ \ln q_j) Cov({\hat q}_i, {\hat q}_j) \\
 & = & (1+\ln q_i)\,(1+ \ln q_j) c^2 v_i \, v_j \, Cov(x_i, x_j).
\end{eqnarray*}
We know that for a multinomial distribution,
\begin{equation}
 Cov(x_i, x_j) = - \frac{p_i\, p_j}{K}.
\end{equation}
Therefore,
\begin{eqnarray}
Cov({\hat q}_i \ln {\hat q}_i, {\hat q}_j \ln {\hat q}_j) & = &
 -  (1+\ln q_i)\,(1+ \ln q_j)\, \frac{q_i\, q_j}{K}.
\end{eqnarray}
So, the variance of ${\hat C}_e$ is
\begin{eqnarray*}
var({\hat C}_e) & = & - \sum_{i=1}^n \sum_{j=1}^n 
\, (1+\ln q_i)\,(1+ \ln q_j)\, \frac{q_i\, q_j}{K} \\
& & + \sum_{i=1}^n (1+\ln q_i)^2 (\frac{c\,v_i\,q_i}{K} -
 \frac{q_i^2}{K})
 + \sum_{i=1}^n (1+\ln q_i)^2 \frac{q_i^2}{K}  \\
& = & - \sum_{i=1}^n \sum_{j=1}^n (1+\ln q_i)\,(1+ \ln q_j)\,
\frac{q_i\, q_j}{K}
 + \sum_{i=1}^n (1+\ln q_i)^2 \frac{c\,v_i\,q_i}{K} \\
& = & - \frac{(1- C_e)^2}{K} 
 + \sum_{i=1}^n (1+\ln q_i)^2 \frac{c\,v_i\,q_i}{K}.
\end{eqnarray*}

\item
From equation~(\ref{eq:ce_mean}), we see that this estimation is
unbiased.
Since $0 < v_i \le 1$ and $q_i$ is non-zero for each $i \in [1,n]$,
and $var({\hat C}_e) \propto \frac{1}{K}$ in
equation~(\ref{eq:ce_variance}), we can have 
\[  \lim_{K \rightarrow \infty} var({\hat C}_e)  = 0. \]
Therefore, we have 
\[  \lim_{K \rightarrow \infty} {\hat C}_e = C_e, \]
and hence ${\hat C}_e$ is a consistent estimator of $C_e$.

\item 
If the population size is large enough that $K \, p_m \ge 5$,
then with confidence of $(1-\alpha) 100 \%$, the estimation error of
the complexity ${\hat C}_e$ is at most $\beta \sigma_{{\hat C}_e}$,
where  $\sigma^2_{{\hat C}_e}=var({\hat C}_e)$. That is,
\begin{eqnarray*}
 | {\hat C}_e - C_e | & \le & \beta\, \sigma_{{\hat C}_e}. \\
\end{eqnarray*}
Replacing $p_i$ with the approximation $x_i$ in
 $\sigma_{{\hat C}_e}$ and considering $0 < v_i \le 1$, we have
\begin{eqnarray*}
\sigma_{{\hat C}_e}^2 & \approx &   - \frac{(1- {\hat C}_e)^2}{K}
 + \sum_{i=1}^n (1+\ln {\hat q}_i)^2 \frac{c\,v_i\,{\hat q}_i}{K} \\
& \le  & -\frac{(1- {\hat C}_e)^2}{K} 
+ \frac{c}{K} \sum_{i=1}^n (1+\ln {\hat q}_i)^2 {\hat q}_i \\ 
 & \le  &  -\frac{(1- {\hat C}_e)^2}{K} 
 + \frac{c}{K} (1 - 2 {\hat C}_e - \ln {\hat q}_{min} 
 {\hat  C}_e) \\ 
& = & -\frac{(1- {\hat C}_e)^2}{K} 
 + \frac{c}{K} (1- (2+\ln {\hat q}_{min}) {\hat C}_e) =
 \delta^2.
\end{eqnarray*}
Since $\sigma_{{\hat C}_e}^2 > 0$, it is always true that 
 $\delta^2 \ge 0$.
We have the third result of this theorem,
\begin{eqnarray*}
   | {\hat C}_e - C_e | & \le & \beta\, \delta.
\end{eqnarray*}

\end{enumerate}

{\Large $\Box  $}

\end{document}